\documentclass[10pt,reqno]{amsart} %pdf-tex
\usepackage{geometry}                
\usepackage{amssymb}
\usepackage{epsfig}
\usepackage{graphics,color}
\usepackage{tikz}

\input{epsf.sty}

\newcommand\web{Sainte-Lagu\"e method}

\newcommand\RED[1]{\color{red} #1}

\newcommand\punkt[1]{\if.#1\else.\spacefactor1000\fi{#1}}

\begingroup
  \divide\count255 by 60
  \multiply\count255 by -60
  \advance\count255 by \time
  \ifnum \count255 < 10 \xdef\klockan{\the\count1.0\the\count255}
  \else\xdef\klockan{\the\count1.\the\count255}\fi
\endgroup

\newcommand{\refS}[1]{Section~\ref{#1}}

\newenvironment{romenumerate}[1][0pt]{% optional argument changes indentation
\addtolength{\leftmargini}{#1}\begin{enumerate}% gives (i), (ii) etc.
 }{\end{enumerate}}

\newtheorem{theorem}{Theorem}[section]

\newtheorem{remark}[theorem]{Remark}

\begin{document}
\title[Dynamic electoral method]
{Dynamic adjustment: \\ an electoral method for relaxed double proportionality}

\author{Svante Linusson} \thanks{Svante Linusson is a Royal Swedish Academy of Sciences Research Fellow supported by a grant from the Knut and Alice Wallenberg Foundation.}
\address{Department of Mathematics, KTH-Royal Institute of Technology, 
  SE-100 44, Stockholm, Sweden.}
\email{linusson@math.kth.se}

\author{Gustav Ryd} 
\address{S:t Eriksgatan 23, SE-11239 Stockholm, Sweden}
\email{gustav.ryd@gmail.com}

\date{\today}  % (typeset \today{} \klockan)} %; revised ...

\begin{abstract}
We describe an electoral system for distributing seats in a parliament. It gives proportionality for the political
parties and close to proportionality for constituencies. The system suggested here 
is a version of the system used in Sweden and other Nordic countries with permanent seats in each 
constituency and adjustment seats to give proportionality on the national level. In the national election
of 2010 the current Swedish system failed to give proportionality between parties. We examine here
one possible cure for this unwanted behavior.
The main difference compared to 
the current Swedish system is that the number of adjustment seats is not fixed, but rather dynamically 
determined to be as low as possible and still insure proportionality between parties.

\end{abstract}

\maketitle

\section{Introduction} \label{S:Intro}

The goal with the current Swedish system is to give proportionality between parties (that pass the $4\%$ threshold) and guarantee each constituency not too few seats, aiming at an approximate proportionality between the 29 constituencies, see \cite {Val}. As described below in \refS{S:Sweden} it can happen that the distribution of seats between the parties deviates from proportionality with some seats for some parties. This 
effect has become more apparent in recent years when the number of parties in the parliament has increased to eight and a parliamentary committee is currently working on changes in the system.

In the present paper a slightly modified system is suggested, see \refS{S:Dynamic}, called the system of dynamic adjustment. It guarantees the proportionality between parties and with this as a side constraint gives reasonable proportionality between constituencies by minimizing the number of adjustment seats. 
The system is based on not having a fixed number of adjustment seats, but rather using 
exactly as few as is needed for that particular election.  In \refS{S:Implement} we describe an efficient way to implement the principle of minimization of adjustment seats and discuss some other practical aspects of an implementation.
In \refS{S:Non-monotone} we scrutinize some theoretical weaknesses of the suggested system, in \refS{S:backtesting} we study how it would have performed in old elections and in 
\refS{S:Simulate} we present simulations done to test the performance. 

Our main conclusion is that this would have worked well in practice in all old elections we have looked at and
would be a possible system that should be seriously considered for future elections. 
One negative aspect of the suggested system 
is a theoretical possibility that extra votes for a party in a constituency might not beneficial for each candidate, but the 
simulations in Section \ref{S:SimNon-monotone} suggest that it is not very likely to occur.  
In \refS{S:Con} we summarize our conclusions. 

\medskip
\noindent
{\it Note added:} On January 10, 2013, a parliamentary committee appointed to suggest changes in the Swedish election act delivered 
their report, \cite{Ni}. The dynamic method discussed here is called Alternativ 5, see page 69 in \cite{Ni}, and was by the 
committee considered to be the second most preferred way to improve the Swedish system.

\medskip
\noindent
{\bf Acknowledgment:} We thank Svante Janson and Jan Lanke for several valuable comments and discussions. The exposition 
has been improved thanks to suggestions from the anonymous referees.

\section{Background} \label{S:Sweden}
Our effort in designing the system here suggested stems from the Swedish electoral system. All official documents in Sweden 
claim the electoral system to be proportional, but this need not be the case. We have done a small, but important, change to the 
Swedish system in order to guarantee proportionality between parties and still preserve a good geographic distribution of seats.

The basis for the Swedish system is \web{} a.k.a. Webster's method  or
the method of odd numbers\footnote{The method seems to be called \web{} in Europe and Webster's method in the US.}. There 
exist equivalent formulations of this method and we will use the version written into Swedish law which we need to describe. The 
seats are given to parties (or constituencies) one at a time, in each stage giving the next seat to the party (constituency) with the 
largest comparison number. The comparison number for party $j$ is obtained as the number of votes for $j$ divided by $2m_j
+1$, where $m_j$ is the number of seats that party $j$ has been given previously in this process.
In Sweden, for historical reasons, the system is (within each constituency) 
used  with the very first divisor being $1.4$ instead of $1$,
which makes the first seat more difficult to get for a small party. This is called the 
{\em modified \web} or modified method of odd numbers.

We find it convenient to use a matrix formulation of the electoral system. We think of the voting as input in the 
form of a matrix $V=(v_{i,j})$ of actual votes, where $v_{i,j}$ is the number of votes party $j$ got in constituency 
$i$. From this we want the electoral system to produce an output matrix $M=(m_{i,j})$ where 
$m_{i,j}$ is the number of seats
party $j$ gets in constituency $i$. The main ingredients in the Swedish system are as follows, for an official 
description with more details see \cite{Val}. The total number of seats in the Swedish parliament (Riksdagen) is 
349.

\begin{romenumerate}
\item{} Before the election, 310 of the seats, called permanent seats, are distributed proportionally to the 29 constituencies (based on entitled voters using Hamilton's method a.k.a the method of largest remainders).
\item{} When the votes have been counted, the permanent seats are distributed to parties within each constituency, i.e. within each row separately, by the modified \web. To participate in the distribution of seats a party must have reached 4\% of the votes nationwide or 12\% of the votes in that particular constituency.
\item{} The next step is to do a theoretical computation of how the 349 seats would have been distributed by
\web{} considering the entire nation as one constituency. This gives for each party $j$ a total number of seats $r_j$ it should have, i.e. the desired column sum of column $j$ in $M$.
\item{} Each party gets extra seats among the 39 adjustment seats so that the total number of seats for each 
party $j$ equals $r_j$. The adjustment seats are distributed within each column with \web{}
in that column, starting from the distribution given by the permanent seats\footnote{Thus, it could very well be that the final outcome in a column is not according to the \web{} if some constituency has been given unproportionally 
many eats previously.}.\\
\item[BUT] it could happen that some party $j$ (or several parties) in Step (ii) received more permanent seats than it should have in total so the column sum exceeds $r_j$. Then party $j$ gets to keep all its permanent seats. The computation in Step (iii) is then redone with the votes and seats of party $j$ removed to find how the remaining seats should be distributed among the remaining parties. Some other party $k$ (or several parties) will then get its number of seats $r_k$ lowered. (It could happen that this needs to be iterated several times.)
\end{romenumerate}

{\bf Example}: In Table \ref{T:2010} is the outcome of the national election in Sweden 2010. In the top row of each cell is the
the number of votes $v_{i,j}$, in the bottom row is  the number of permanent seats received to the left and the number of
adjustment seats to the right (and in red). The rightmost column shows the number of entitled voters for each constituency.
We see that the number of permanent seats for (M) and (S), second to last row 107 and 112 in step (ii) 
are more than the number of total seats they should have, see last row 106 and 109 respectively. 
Thus, according to the BUT, step number (iii) has to be redone which 
gives the new totals for the other parties, shown to the right in the second line of second to last row. This determines the number of adjustment seats for each party to be distributed in its column.
In the second line of the last row is the amount of disproportionality for each party.
The parties in the Swedish parliament are given by the standard ordering used in Sweden, first the four center-right parties (that currently form a minority government): (M) Moderaterna - conservatives, (C) Centern - agrarian/liberals, (FP) Folkpartiet - liberals, (KD) Kristdemokraterna - christian democrats, followed by (S) Socialdemokraterna - social democrats, (V) V\"anstern - left, former communists, (MP) Milj\"opartiet - greens, and last the newest party in the Swedish parliament: (SD) Sverigedemokraterna - Swedish democrats, nationalistic/anti-immigration party.

\begin{table}[h]
\renewcommand{\arraystretch}{0.68}
\begin{tabular}{| l || lr | lr | lr  | lr  | lr  | lr  | lr  | lr  || lr || }
\hline
{\tiny   Constituency $\backslash$ Party }&\multicolumn{2}{|c|}{\tiny (M) }&\multicolumn{2}{|c|}{\tiny (C) }&\multicolumn{2}{|c|}{\tiny (FP)}&\multicolumn{2}{|c|}{\tiny (KD) }&
\multicolumn{2}{|c|} {\tiny (S)}&\multicolumn{2}{|c|}{\tiny (V)}&\multicolumn{2}{|c|} {\tiny (MP)}&\multicolumn{2}{|c||} {\tiny (SD) }&\multicolumn{2}{|c||} {\tiny \hspace{-6pt}Entitled voters } \\
\hline
\hline
{\tiny    Stockholms stad }& \multicolumn{2}{|r|}{\tiny \hspace{-6pt} 183 421}&\multicolumn{2}{|r|}{\tiny  33 895}&
\multicolumn{2}{|r|}{\tiny   45 939}&\multicolumn{2}{|r|}{\tiny  28 244}&
\multicolumn{2}{|r|}{\tiny  \hspace{-6pt} 111 688}&\multicolumn{2}{|r|}{\tiny 39 565}&\multicolumn{2}{|r|} {\tiny 65 351}&
\multicolumn{2}{|r||}{\tiny 16 950}&\multicolumn{2}{|r||} {\tiny \hspace{-6pt}634 464} \\
{\tiny }&{\tiny 10}&{\tiny }&{\tiny  2}&{\tiny }&{\tiny 2}&\RED {\tiny +1}&{\tiny 2 }&{\tiny }&
{\tiny 6}&{\tiny }&{\tiny 2}&{\tiny }& {\tiny 3}&{\tiny }&{\tiny 1}&{\tiny }&{\tiny 29}&\RED{\tiny +1} \\
\hline

{\tiny  Stockholms l\"an  }& \multicolumn{2}{|r|}{\tiny  \hspace{-6pt} 286 249}&\multicolumn{2}{|r|}{\tiny 41 369 }&
\multicolumn{2}{|r|}{\tiny 59 461 }&\multicolumn{2}{|r|}{\tiny 44 880 }&
\multicolumn{2}{|r|}{\tiny  \hspace{-6pt} 159 222}&\multicolumn{2}{|r|}{\tiny 31 617}&\multicolumn{2}{|r|} {\tiny 53 788}&
\multicolumn{2}{|r||}{\tiny 29 886}&\multicolumn{2}{|r||} {\tiny \hspace{-6pt}850 629 } \\
{\tiny }&{\tiny 15}&\RED{\tiny }&{\tiny  2}&\RED{\tiny }&{\tiny 3}&\RED {\tiny }&{\tiny  2}&\RED{\tiny +1}&
{\tiny 8}&\RED{\tiny }&{\tiny 2}&\RED {\tiny }& {\tiny 3}&\RED {\tiny }&{\tiny 2}&\RED {\tiny }&{\tiny 38}&\RED {\tiny +1} \\
\hline

{\tiny   Uppsala l\"an }& \multicolumn{2}{|r|}{\tiny 64 750}&\multicolumn{2}{|r|}{\tiny 17 838 }&
\multicolumn{2}{|r|}{\tiny 16 878 }&\multicolumn{2}{|r|}{\tiny  12 265}&
\multicolumn{2}{|r|}{\tiny 58 862  }&\multicolumn{2}{|r|}{\tiny 11 845}&\multicolumn{2}{|r|} {\tiny 18 993}&
\multicolumn{2}{|r||}{\tiny 10 003}&\multicolumn{2}{|r||} {\tiny \hspace{-6pt}253 765 } \\
{\tiny }&{\tiny 4}&\RED{\tiny }&{\tiny 1 }&\RED{\tiny }&{\tiny 1}&\RED {\tiny }&{\tiny  1}&\RED{\tiny }&
{\tiny 3}&\RED{\tiny }&{\tiny }&\RED {\tiny +1}& {\tiny 1}&\RED {\tiny }&{\tiny }&\RED {\tiny +1}&{\tiny 11}&\RED {\tiny +2} \\
\hline

{\tiny  S\"odermanlands l\"an  }& \multicolumn{2}{|r|}{\tiny 47 889}&\multicolumn{2}{|r|}{\tiny  9 850}&
\multicolumn{2}{|r|}{\tiny 11 299 }&\multicolumn{2}{|r|}{\tiny 8 095 }&
\multicolumn{2}{|r|}{\tiny  59 463 }&\multicolumn{2}{|r|}{\tiny 8 637}&\multicolumn{2}{|r|} {\tiny 13 065 }&
\multicolumn{2}{|r||}{\tiny 11 370}&\multicolumn{2}{|r||} {\tiny \hspace{-6pt}204 779} \\
{\tiny }&{\tiny 3}&\RED{\tiny }&{\tiny  }&\RED{\tiny +1}&{\tiny }&\RED {\tiny +1}&{\tiny  }&\RED{\tiny }&
{\tiny 4}&\RED{\tiny }&{\tiny }&\RED {\tiny }& {\tiny 1}&\RED {\tiny }&{\tiny 1}&\RED {\tiny }&{\tiny 9}&\RED {\tiny +2} \\
\hline

{\tiny  \"Osterg\"otlands l\"an  }& \multicolumn{2}{|r|}{\tiny 80 141}&\multicolumn{2}{|r|}{\tiny 17 561 }&
\multicolumn{2}{|r|}{\tiny  19 017}&\multicolumn{2}{|r|}{\tiny 16 357}&
\multicolumn{2}{|r|}{\tiny 92 164}&\multicolumn{2}{|r|}{\tiny 14 242}&\multicolumn{2}{|r|} {\tiny 21 225}&
\multicolumn{2}{|r||}{\tiny 14 862}&\multicolumn{2}{|r||} {\tiny \hspace{-6pt}330 010} \\
{\tiny }&{\tiny 4}&\RED{\tiny }&{\tiny 1}&\RED{\tiny }&{\tiny 1}&\RED {\tiny }&{\tiny 1}&\RED{\tiny }&
{\tiny 5}&\RED{\tiny }&{\tiny }&\RED {\tiny +1}& {\tiny 1}&\RED {\tiny }&{\tiny 1}&\RED {\tiny }&{\tiny 14}&\RED {\tiny +1} \\
\hline

{\tiny   J\"onk\"opings l\"an}& \multicolumn{2}{|r|}{\tiny 57 901}&\multicolumn{2}{|r|}{\tiny  16 859}&
\multicolumn{2}{|r|}{\tiny  12 134}&\multicolumn{2}{|r|}{\tiny  27 822}&
\multicolumn{2}{|r|}{\tiny  66 316}&\multicolumn{2}{|r|}{\tiny 8 775}&\multicolumn{2}{|r|} {\tiny 11 438}&
\multicolumn{2}{|r||}{\tiny 13 888}&\multicolumn{2}{|r||} {\tiny \hspace{-6pt}256 538 } \\
{\tiny }&{\tiny 3}&\RED{\tiny }&{\tiny 1}&\RED{\tiny }&{\tiny }&\RED {\tiny +1}&{\tiny 2}&\RED{\tiny }&
{\tiny 4}&\RED{\tiny }&{\tiny }&\RED {\tiny }& {\tiny }&\RED {\tiny +1}&{\tiny 1}&\RED {\tiny }&{\tiny 11}&\RED {\tiny +2} \\
\hline

{\tiny Kronobergs l\"an}& \multicolumn{2}{|r|}{\tiny 34 762}&\multicolumn{2}{|r|}{\tiny 11 559}&
\multicolumn{2}{|r|}{\tiny 6 667}&\multicolumn{2}{|r|}{\tiny 7 111}&
\multicolumn{2}{|r|}{\tiny  35 555}&\multicolumn{2}{|r|}{\tiny 5 380}&\multicolumn{2}{|r|} {\tiny 7 044}&
\multicolumn{2}{|r||}{\tiny 7 424}&\multicolumn{2}{|r||} {\tiny \hspace{-6pt}138 781} \\
{\tiny }&{\tiny 2}&\RED{\tiny }&{\tiny 1}&\RED{\tiny }&{\tiny }&\RED {\tiny }&{\tiny  }&\RED{\tiny }&
{\tiny 3}&\RED{\tiny }&{\tiny }&\RED {\tiny }& {\tiny }&\RED {\tiny }&{\tiny }&\RED {\tiny }&{\tiny 6}&\RED {\tiny } \\
\hline

{\tiny Kalmar l\"an}& \multicolumn{2}{|r|}{\tiny 41 631}&\multicolumn{2}{|r|}{\tiny 13 829}&
\multicolumn{2}{|r|}{\tiny 7 847}&\multicolumn{2}{|r|}{\tiny 9 341}&
\multicolumn{2}{|r|}{\tiny 55 116}&\multicolumn{2}{|r|}{\tiny 7 679}&\multicolumn{2}{|r|} {\tiny 8 713}&
\multicolumn{2}{|r||}{\tiny 8 964}&\multicolumn{2}{|r||} {\tiny \hspace{-6pt}184 737} \\
{\tiny }&{\tiny 3}&\RED{\tiny }&{\tiny 1}&\RED{\tiny }&{\tiny }&\RED {\tiny }&{\tiny  }&\RED{\tiny +1}&
{\tiny 4}&\RED{\tiny }&{\tiny }&\RED {\tiny }& {\tiny }&\RED {\tiny }&{\tiny }&\RED {\tiny }&{\tiny 8}&\RED {\tiny +1} \\
\hline

{\tiny Gotlands l\"an}& \multicolumn{2}{|r|}{\tiny 9 731}&\multicolumn{2}{|r|}{\tiny 5 657}&
\multicolumn{2}{|r|}{\tiny 1 785}&\multicolumn{2}{|r|}{\tiny 1 128}&
\multicolumn{2}{|r|}{\tiny  12 855}&\multicolumn{2}{|r|}{\tiny 2 342}&\multicolumn{2}{|r|} {\tiny 3 259}&
\multicolumn{2}{|r||}{\tiny 1 225}&\multicolumn{2}{|r||} {\tiny \hspace{-6pt}46 237} \\
{\tiny }&{\tiny 1}&\RED{\tiny }&{\tiny  }&\RED{\tiny }&{\tiny }&\RED {\tiny }&{\tiny  }&\RED{\tiny }&
{\tiny 1}&\RED{\tiny }&{\tiny }&\RED {\tiny }& {\tiny }&\RED {\tiny }&{\tiny }&\RED {\tiny }&{\tiny 2}&\RED {\tiny } \\
\hline

{\tiny Blekinge l\"an}& \multicolumn{2}{|r|}{\tiny 27 387}&\multicolumn{2}{|r|}{\tiny 5 771}&
\multicolumn{2}{|r|}{\tiny 5 431}&\multicolumn{2}{|r|}{\tiny 3 973}&
\multicolumn{2}{|r|}{\tiny 36 520}&\multicolumn{2}{|r|}{\tiny 5 075}&\multicolumn{2}{|r|} {\tiny 5 289}&
\multicolumn{2}{|r||}{\tiny 9 830}&\multicolumn{2}{|r||} {\tiny \hspace{-6pt}118 279} \\
{\tiny }&{\tiny 2}&\RED{\tiny }&{\tiny  }&\RED{\tiny }&{\tiny }&\RED {\tiny }&{\tiny  }&\RED{\tiny }&
{\tiny 3}&\RED{\tiny }&{\tiny }&\RED {\tiny }& {\tiny }&\RED {\tiny }&{\tiny }&\RED {\tiny +1}&{\tiny 5}&\RED {\tiny +1} \\
\hline

{\tiny Malm\"o kommun}& \multicolumn{2}{|r|}{\tiny 55 160}&\multicolumn{2}{|r|}{\tiny 4 795}&
\multicolumn{2}{|r|}{\tiny 11 768}&\multicolumn{2}{|r|}{\tiny 5 274}&
\multicolumn{2}{|r|}{\tiny 48 450}&\multicolumn{2}{|r|}{\tiny 10 118}&\multicolumn{2}{|r|} {\tiny 14 861}&
\multicolumn{2}{|r||}{\tiny 13 256}&\multicolumn{2}{|r||} {\tiny \hspace{-6pt}214 326} \\
{\tiny }&{\tiny 3}&\RED{\tiny }&{\tiny  }&\RED{\tiny }&{\tiny 1}&\RED {\tiny }&{\tiny  }&\RED{\tiny }&
{\tiny 3}&\RED{\tiny }&{\tiny }&\RED {\tiny +1}& {\tiny 1}&\RED {\tiny }&{\tiny 1}&\RED {\tiny }&{\tiny 9}&\RED {\tiny +1} \\
\hline

{\tiny Sk{\aa}ne l\"ans }& \multicolumn{2}{|r|}{\tiny 58 628}&\multicolumn{2}{|r|}{\tiny 8 164}&
\multicolumn{2}{|r|}{\tiny 13 967}&\multicolumn{2}{|r|}{\tiny  6 989}&
\multicolumn{2}{|r|}{\tiny 49 900}&\multicolumn{2}{|r|}{\tiny 5 847}&\multicolumn{2}{|r|} {\tiny 9 869}&
\multicolumn{2}{|r||}{\tiny 17 448}&\multicolumn{2}{|r||} {\tiny \hspace{-6pt}213 580} \\
{\tiny v\"astra}&{\tiny 4}&\RED{\tiny }&{\tiny  }&\RED{\tiny }&{\tiny 1}&\RED {\tiny }&{\tiny  }&\RED{\tiny }&
{\tiny 3}&\RED{\tiny }&{\tiny }&\RED {\tiny }& {\tiny }&\RED {\tiny +1}&{\tiny 1}&\RED {\tiny }&{\tiny 9}&\RED {\tiny +1} \\
\hline

{\tiny Sk{\aa}ne l\"ans }& \multicolumn{2}{|r|}{\tiny 87 893}&\multicolumn{2}{|r|}{\tiny 12 717}&
\multicolumn{2}{|r|}{\tiny 19 622}&\multicolumn{2}{|r|}{\tiny 9 916}&
\multicolumn{2}{|r|}{\tiny 50 557}&\multicolumn{2}{|r|}{\tiny 7 597}&\multicolumn{2}{|r|} {\tiny 16 176}&
\multicolumn{2}{|r||}{\tiny 19 923}&\multicolumn{2}{|r||} {\tiny \hspace{-6pt}267 562} \\
{\tiny s\"odra}&{\tiny 5}&\RED{\tiny }&{\tiny 1}&\RED{\tiny }&{\tiny 1}&\RED {\tiny }&{\tiny  }&\RED{\tiny +1}&
{\tiny 3}&\RED{\tiny }&{\tiny }&\RED {\tiny }& {\tiny 1}&\RED {\tiny }&{\tiny 1}&\RED {\tiny }&{\tiny 12}&\RED {\tiny +1} \\
\hline

{\tiny   Sk{\aa}ne l\"ans}& \multicolumn{2}{|r|}{\tiny 60 930}&\multicolumn{2}{|r|}{\tiny 12 871}&
\multicolumn{2}{|r|}{\tiny 12 677}&\multicolumn{2}{|r|}{\tiny 9 420}&
\multicolumn{2}{|r|}{\tiny 54 529}&\multicolumn{2}{|r|}{\tiny 6 113}&\multicolumn{2}{|r|} {\tiny 10 195}&
\multicolumn{2}{|r||}{\tiny 21 312}&\multicolumn{2}{|r||} {\tiny \hspace{-6pt}232 273} \\
{\tiny norra och \"ostra}&{\tiny 4}&\RED{\tiny }&{\tiny 1}&\RED{\tiny }&{\tiny 1}&\RED {\tiny }&{\tiny  }&\RED{\tiny +1}&
{\tiny 3}&\RED{\tiny }&{\tiny }&\RED {\tiny }& {\tiny }&\RED {\tiny +1}&{\tiny 1}&\RED {\tiny }&{\tiny 10}&\RED {\tiny +2} \\
\hline

{\tiny Hallands l\"an}& \multicolumn{2}{|r|}{\tiny 67 878}&\multicolumn{2}{|r|}{\tiny 17 178}&
\multicolumn{2}{|r|}{\tiny 15 286}&\multicolumn{2}{|r|}{\tiny 10 994}&
\multicolumn{2}{|r|}{\tiny 52 319}&\multicolumn{2}{|r|}{\tiny 6 904}&\multicolumn{2}{|r|} {\tiny 11 568}&
\multicolumn{2}{|r||}{\tiny 10 507}&\multicolumn{2}{|r||} {\tiny \hspace{-6pt}229 891} \\
{\tiny }&{\tiny 4}&\RED{\tiny }&{\tiny 1}&\RED{\tiny }&{\tiny 1}&\RED {\tiny }&{\tiny  }&\RED{\tiny +1}&
{\tiny 3}&\RED{\tiny }&{\tiny }&\RED {\tiny }& {\tiny 1}&\RED {\tiny }&{\tiny }&\RED {\tiny +1}&{\tiny 10}&\RED {\tiny +2} \\
\hline

{\tiny G\"oteborgs kommun}& \multicolumn{2}{|r|}{\tiny 96 981}&\multicolumn{2}{|r|}{\tiny 12 183}&
\multicolumn{2}{|r|}{\tiny 26 829}&\multicolumn{2}{|r|}{\tiny 19 484}&
\multicolumn{2}{|r|}{\tiny 80 543}&\multicolumn{2}{|r|}{\tiny 27 246}&\multicolumn{2}{|r|} {\tiny 34 205}&
\multicolumn{2}{|r||}{\tiny 15 608}&\multicolumn{2}{|r||} {\tiny \hspace{-6pt}389 821} \\
{\tiny }&{\tiny 5}&\RED{\tiny }&{\tiny  }&\RED{\tiny +1}&{\tiny 1}&\RED {\tiny }&{\tiny 1}&\RED{\tiny }&
{\tiny 5}&\RED{\tiny }&{\tiny 2}&\RED {\tiny }& {\tiny 2}&\RED {\tiny }&{\tiny 1}&\RED {\tiny }&{\tiny 17}&\RED {\tiny +1} \\
\hline

{\tiny V G\"otalands l\"an}& \multicolumn{2}{|r|}{\tiny 73 853}&\multicolumn{2}{|r|}{\tiny 13 563}&
\multicolumn{2}{|r|}{\tiny 20 194}&\multicolumn{2}{|r|}{\tiny 16 525}&
\multicolumn{2}{|r|}{\tiny 59 477}&\multicolumn{2}{|r|}{\tiny 10 506}&\multicolumn{2}{|r|} {\tiny 15 794}&
\multicolumn{2}{|r||}{\tiny 12 504}&\multicolumn{2}{|r||} {\tiny \hspace{-6pt}264 666} \\
{\tiny v\"astra}&{\tiny 4}&\RED{\tiny }&{\tiny 1}&\RED{\tiny }&{\tiny 1}&\RED {\tiny }&{\tiny 1}&\RED{\tiny }&
{\tiny 3}&\RED{\tiny }&{\tiny }&\RED {\tiny +1}& {\tiny 1}&\RED {\tiny }&{\tiny 1}&\RED {\tiny }&{\tiny 12}&\RED {\tiny +1} \\
\hline

{\tiny V G\"otalands l\"an}& \multicolumn{2}{|r|}{\tiny 46 582}&\multicolumn{2}{|r|}{\tiny 11 449}&
\multicolumn{2}{|r|}{\tiny 13 393}&\multicolumn{2}{|r|}{\tiny 11 092}&
\multicolumn{2}{|r|}{\tiny 56 060}&\multicolumn{2}{|r|}{\tiny 9 907}&\multicolumn{2}{|r|} {\tiny 12 003}&
\multicolumn{2}{|r||}{\tiny 10 513}&\multicolumn{2}{|r||} {\tiny \hspace{-6pt}205 328} \\
{\tiny norra}&{\tiny 3}&\RED{\tiny }&{\tiny 1}&\RED{\tiny }&{\tiny 1}&\RED {\tiny }&{\tiny  }&\RED{\tiny +1}&
{\tiny 3}&\RED{\tiny }&{\tiny }&\RED {\tiny +1}& {\tiny 1}&\RED {\tiny }&{\tiny }&\RED {\tiny +1}&{\tiny 9}&\RED {\tiny +3} \\
\hline

{\tiny V G\"otalands l\"an}& \multicolumn{2}{|r|}{\tiny 34 334}&\multicolumn{2}{|r|}{\tiny 9 273}&
\multicolumn{2}{|r|}{\tiny 8 883}&\multicolumn{2}{|r|}{\tiny 7 745}&
\multicolumn{2}{|r|}{\tiny 37 817}&\multicolumn{2}{|r|}{\tiny 6 136}&\multicolumn{2}{|r|} {\tiny 7 315}&
\multicolumn{2}{|r||}{\tiny 8 350}&\multicolumn{2}{|r||} {\tiny \hspace{-6pt}144 186} \\
{\tiny s\"odra}&{\tiny 3}&\RED{\tiny }&{\tiny  }&\RED{\tiny }&{\tiny }&\RED {\tiny }&{\tiny  }&\RED{\tiny }&
{\tiny 3}&\RED{\tiny }&{\tiny }&\RED {\tiny }& {\tiny }&\RED {\tiny }&{\tiny }&\RED {\tiny }&{\tiny 6}&\RED {\tiny } \\
\hline

{\tiny V G\"otalands l\"an}& \multicolumn{2}{|r|}{\tiny 47 049}&\multicolumn{2}{|r|}{\tiny 13 914}&
\multicolumn{2}{|r|}{\tiny 10 387}&\multicolumn{2}{|r|}{\tiny 11 092}&
\multicolumn{2}{|r|}{\tiny 57 095}&\multicolumn{2}{|r|}{\tiny 8 223}&\multicolumn{2}{|r|} {\tiny 9 440}&
\multicolumn{2}{|r||}{\tiny 9 725}&\multicolumn{2}{|r||} {\tiny \hspace{-6pt}200 322} \\
{\tiny \"ostra}&{\tiny 3}&\RED{\tiny }&{\tiny 1}&\RED{\tiny }&{\tiny }&\RED {\tiny +1}&{\tiny  1}&\RED{\tiny }&
{\tiny 4}&\RED{\tiny }&{\tiny }&\RED {\tiny }& {\tiny }&\RED {\tiny }&{\tiny }&\RED {\tiny }&{\tiny 9}&\RED {\tiny +1} \\
\hline

{\tiny V\"armlands l\"an}& \multicolumn{2}{|r|}{\tiny 45 578}&\multicolumn{2}{|r|}{\tiny 13 379}&
\multicolumn{2}{|r|}{\tiny 10 652}&\multicolumn{2}{|r|}{\tiny 8 312}&
\multicolumn{2}{|r|}{\tiny 68 520}&\multicolumn{2}{|r|}{\tiny 10 231}&\multicolumn{2}{|r|} {\tiny 9 997}&
\multicolumn{2}{|r||}{\tiny 8 502}&\multicolumn{2}{|r||} {\tiny \hspace{-6pt}213 239} \\
{\tiny }&{\tiny 3}&\RED{\tiny }&{\tiny 1}&\RED{\tiny }&{\tiny }&\RED {\tiny +1}&{\tiny  }&\RED{\tiny }&
{\tiny 5}&\RED{\tiny }&{\tiny }&\RED {\tiny +1}& {\tiny }&\RED {\tiny +1}&{\tiny }&\RED {\tiny }&{\tiny 9}&\RED {\tiny +3} \\
\hline

{\tiny \"Orebro l\"an}& \multicolumn{2}{|r|}{\tiny 43 791}&\multicolumn{2}{|r|}{\tiny 9 807}&
\multicolumn{2}{|r|}{\tiny 11 415}&\multicolumn{2}{|r|}{\tiny 11 235}&
\multicolumn{2}{|r|}{\tiny 70 818}&\multicolumn{2}{|r|}{\tiny 10 311}&\multicolumn{2}{|r|} {\tiny 11 846}&
\multicolumn{2}{|r||}{\tiny 11 136}&\multicolumn{2}{|r||} {\tiny \hspace{-6pt}215 772} \\
{\tiny }&{\tiny 3}&\RED{\tiny }&{\tiny  }&\RED{\tiny }&{\tiny 1}&\RED {\tiny }&{\tiny  }&\RED{\tiny +1}&
{\tiny 4}&\RED{\tiny }&{\tiny }&\RED {\tiny +1}& {\tiny 1}&\RED {\tiny }&{\tiny }&\RED {\tiny +1}&{\tiny 9}&\RED {\tiny +3} \\
\hline

{\tiny V\"astmanlands l\"an}& \multicolumn{2}{|r|}{\tiny 43 462}&\multicolumn{2}{|r|}{\tiny 8 266}&
\multicolumn{2}{|r|}{\tiny 12 016}&\multicolumn{2}{|r|}{\tiny 7 406}&
\multicolumn{2}{|r|}{\tiny 58 222}&\multicolumn{2}{|r|}{\tiny 9 154}&\multicolumn{2}{|r|} {\tiny 9 459}&
\multicolumn{2}{|r||}{\tiny 9 992}&\multicolumn{2}{|r||} {\tiny \hspace{-6pt}192 258} \\
{\tiny }&{\tiny 3}&\RED{\tiny }&{\tiny  }&\RED{\tiny }&{\tiny 1}&\RED {\tiny }&{\tiny  }&\RED{\tiny }&
{\tiny 4}&\RED{\tiny }&{\tiny }&\RED {\tiny +1}& {\tiny }&\RED {\tiny +1}&{\tiny }&\RED {\tiny +1}&{\tiny 8}&\RED {\tiny +3} \\
\hline

{\tiny Dalarnas l\"an}& \multicolumn{2}{|r|}{\tiny 44 997}&\multicolumn{2}{|r|}{\tiny 14 086}&
\multicolumn{2}{|r|}{\tiny 8 747}&\multicolumn{2}{|r|}{\tiny 7 925}&
\multicolumn{2}{|r|}{\tiny  67 139}&\multicolumn{2}{|r|}{\tiny 10 533}&\multicolumn{2}{|r|} {\tiny 10 652}&
\multicolumn{2}{|r||}{\tiny 12 470}&\multicolumn{2}{|r||} {\tiny \hspace{-6pt}217 072} \\
{\tiny }&{\tiny 3}&\RED{\tiny }&{\tiny  1}&\RED{\tiny }&{\tiny }&\RED {\tiny }&{\tiny  }&\RED{\tiny }&
{\tiny 4}&\RED{\tiny }&{\tiny }&\RED {\tiny +1}& {\tiny 1}&\RED {\tiny }&{\tiny 1}&\RED {\tiny }&{\tiny 10}&\RED {\tiny +1} \\
\hline

{\tiny G\"avleborgs l\"an}& \multicolumn{2}{|r|}{\tiny 41 009}&\multicolumn{2}{|r|}{\tiny 12 982}&
\multicolumn{2}{|r|}{\tiny 9 444}&\multicolumn{2}{|r|}{\tiny 7 235}&
\multicolumn{2}{|r|}{\tiny 67 893}&\multicolumn{2}{|r|}{\tiny 12 814}&\multicolumn{2}{|r|} {\tiny 10 918}&
\multicolumn{2}{|r||}{\tiny 12 616}&\multicolumn{2}{|r||} {\tiny \hspace{-6pt} 217 152} \\
{\tiny }&{\tiny 3}&\RED{\tiny }&{\tiny  1}&\RED{\tiny }&{\tiny }&\RED {\tiny +1}&{\tiny  }&\RED{\tiny }&
{\tiny 4}&\RED{\tiny }&{\tiny 1}&\RED {\tiny }& {\tiny }&\RED {\tiny +1}&{\tiny 1}&\RED {\tiny }&{\tiny 10}&\RED {\tiny +2} \\
\hline

{\tiny  V\"asternorrlands l\"an}& \multicolumn{2}{|r|}{\tiny 34 550}&\multicolumn{2}{|r|}{\tiny 11 185}&
\multicolumn{2}{|r|}{\tiny 8 253}&\multicolumn{2}{|r|}{\tiny 6 983}&
\multicolumn{2}{|r|}{\tiny 70 341}&\multicolumn{2}{|r|}{\tiny 9 642}&\multicolumn{2}{|r|} {\tiny 8 757}&
\multicolumn{2}{|r||}{\tiny 7 264}&\multicolumn{2}{|r||} {\tiny \hspace{-6pt}191 150} \\
{\tiny }&{\tiny 2}&\RED{\tiny }&{\tiny 1 }&\RED{\tiny }&{\tiny }&\RED {\tiny }&{\tiny  }&\RED{\tiny }&
{\tiny 5}&\RED{\tiny }&{\tiny }&\RED {\tiny +1}& {\tiny }&\RED {\tiny }&{\tiny }&\RED {\tiny }&{\tiny 8}&\RED {\tiny +1} \\
\hline

{\tiny J\"amtlands l\"an}& \multicolumn{2}{|r|}{\tiny 18 193}&\multicolumn{2}{|r|}{\tiny 10 487}&
\multicolumn{2}{|r|}{\tiny 3 155}&\multicolumn{2}{|r|}{\tiny 2 340}&
\multicolumn{2}{|r|}{\tiny  33 013}&\multicolumn{2}{|r|}{\tiny 5 340}&\multicolumn{2}{|r|} {\tiny 5 339}&
\multicolumn{2}{|r||}{\tiny 3 122}&\multicolumn{2}{|r||} {\tiny \hspace{-6pt}100 144} \\
{\tiny }&{\tiny 1}&\RED{\tiny }&{\tiny 1}&\RED{\tiny }&{\tiny }&\RED {\tiny }&{\tiny  }&\RED{\tiny }&
{\tiny 2}&\RED{\tiny }&{\tiny }&\RED {\tiny }& {\tiny }&\RED {\tiny }&{\tiny }&\RED {\tiny }&{\tiny 4}&\RED {\tiny } \\
\hline

{\tiny V\"asterbottens l\"an}& \multicolumn{2}{|r|}{\tiny 30 184}&\multicolumn{2}{|r|}{\tiny 12 699}&
\multicolumn{2}{|r|}{\tiny 10 296}&\multicolumn{2}{|r|}{\tiny 9 125 }&
\multicolumn{2}{|r|}{\tiny  72 008}&\multicolumn{2}{|r|}{\tiny 17 034}&\multicolumn{2}{|r|} {\tiny 12 246}&
\multicolumn{2}{|r||}{\tiny 4 651}&\multicolumn{2}{|r||} {\tiny \hspace{-6pt} 201 902} \\
{\tiny }&{\tiny 2}&\RED{\tiny }&{\tiny 1}&\RED{\tiny }&{\tiny }&\RED {\tiny +1}&{\tiny  }&\RED{\tiny +1}&
{\tiny 4}&\RED{\tiny }&{\tiny 1}&\RED {\tiny }& {\tiny 1}&\RED {\tiny }&{\tiny }&\RED {\tiny }&{\tiny 9}&\RED {\tiny +2} \\
\hline

{\tiny  Norrbottens l\"an}& \multicolumn{2}{|r|}{\tiny 26 852}&\multicolumn{2}{|r|}{\tiny 7 618}&
\multicolumn{2}{|r|}{\tiny 7 082}&\multicolumn{2}{|r|}{\tiny 5 388}&
\multicolumn{2}{|r|}{\tiny 85 035}&\multicolumn{2}{|r|}{\tiny 15 240}&\multicolumn{2}{|r|} {\tiny 8 630}&
\multicolumn{2}{|r||}{\tiny 6 309}&\multicolumn{2}{|r||} {\tiny \hspace{-6pt} 194 788} \\
{\tiny }&{\tiny 2}&\RED{\tiny }&{\tiny  }&\RED{\tiny }&{\tiny }&\RED {\tiny }&{\tiny  }&\RED{\tiny }&
{\tiny 6}&\RED{\tiny }&{\tiny 1}&\RED {\tiny }& {\tiny }&\RED {\tiny }&{\tiny }&\RED {\tiny }&{\tiny 9}&\RED {\tiny } \\
\hline\hline

{\tiny  Total  }& \multicolumn{2}{|r|}{\tiny \hspace{-4pt}1 791 766}&\multicolumn{2}{|r|}{\tiny \hspace{-6pt}390 804}&
\multicolumn{2}{|r|}{\tiny \hspace{-6pt}420 524}&\multicolumn{2}{|r|}{\tiny  \hspace{-6pt}333 696}&
\multicolumn{2}{|r|}{\tiny  \hspace{-4pt}1 827 497}&\multicolumn{2}{|r|}{\tiny \hspace{-6pt}334 053}&\multicolumn{2}{|r|} {\tiny\hspace{-6pt} 437 435}&
\multicolumn{2}{|r||}{\tiny \hspace{-6pt}339 610}&\multicolumn{2}{|r||} {\tiny \hspace{-4pt}7 123 651} \\
{\tiny }&{\tiny \hspace{-4pt}107}&\RED{\tiny }&\multicolumn{2}{|r|}{{\tiny \hspace{-4pt}21}{\RED\tiny +2}{\tiny =23}}&
\multicolumn{2}{|r|}{{\tiny \hspace{-4pt}17}{\RED\tiny +7}{\tiny =24}}&\multicolumn{2}{|r|}{{\tiny\hspace{-4pt}11}{\RED\tiny +8}{\tiny =19}}&
{\tiny \hspace{-4pt}112}&{\RED\tiny }&\multicolumn{2}{|r|}{{\tiny \hspace{-4pt}9}{\RED\tiny +10}{\tiny =19}}&
\multicolumn{2}{|r|}{{\tiny \hspace{-4pt}19}{\RED\tiny +6}{\tiny =25}}&\multicolumn{2}{|r||}{{\tiny \hspace{-4pt}14}{\RED\tiny +6}{\tiny =20}}&
\multicolumn{2}{|r||}{{\tiny \hspace{-4pt}310}{\RED\tiny +39}{\tiny =349}} \\
\hline

{\tiny Proportionality}& \multicolumn{2}{|r|}{\tiny 106}&\multicolumn{2}{|r|}{\tiny 23 }&
\multicolumn{2}{|r|}{\tiny 25}&\multicolumn{2}{|r|}{\tiny  20}&
\multicolumn{2}{|r|}{\tiny 109 }&\multicolumn{2}{|r|}{\tiny 20}&\multicolumn{2}{|r|} {\tiny 26}&
\multicolumn{2}{|r||}{\tiny 20}&\multicolumn{2}{|r||} {\tiny \hspace{-6pt} } \\
{\tiny Difference }&{\tiny 1}&\RED{\tiny }&{\tiny  }&\RED{\tiny }&{\tiny -1}&\RED {\tiny }&{\tiny  -1}&\RED{\tiny }&
{\tiny 3}&\RED{\tiny }&{\tiny -1}&\RED {\tiny }& {\tiny -1}&\RED {\tiny }&{\tiny }&\RED {\tiny }&{\tiny }&\RED {\tiny } \\
\hline
\end{tabular}
\vskip1mm
\caption{Outcome of Swedish national election 2010.} \label{T:2010}
\end{table}

\medskip
The problem of achieving double proportionality is common for many countries and has found different 
practical solutions in different countries, see e.g. \cite{EU}. Several countries still have only what we here call 
permanent seats, which guarantees good proportionality geographically, but can lead to poor proportionality 
between parties. In other countries the distribution of seats between parties is dominant over geographic
distribution, e.g. Israel and the Netherlands.  It is our perception that in Sweden, although geographic 
spread of members of parliament is important,
the fairness achieved by good proportionality between parties is of highest importance. This was also 
explicitly the
highest priority by the committee that once formed the Swedish system in 1967, see \cite{Sydow}.

\subsection{Elections 2010} In the national Swedish election to Riksdagen in September 2010 the BUT happened to two parties that got too many seats. 
Moderaterna (M) got 1 seat too many and Socialdemokraterna (S) got 3 seats too many, see 
Table 1. 
Since the election was very close it almost affected the majority situation in the parliament and a small fraction\footnote{ Less than 2\%.} of votes, the absentee ballots, 
are counted for the first time on the Wednesday three days after the election this became very apparent and 
attracted large public attention. 
All official descriptions of the Swedish electoral system, e.g. \cite{Val}, start with saying that it is proportional, which is also clearly the intention, the possibility of the BUT came as a surprise to many Swedes. After the current electoral system was introduced in 1970, the BUT has happened once before in the national elections but it was not close to influencing the power balance in the parliament so it was not widely noticed.
 
In fact, by the preliminary results 2010 (without absentee ballots) the Social democrats (S) had received 4 seats too many, but after counting absentee ballots (completed three days after the election) they lost a permanent seat to the Greens (MP) in the constituency of Dalarna, because (MP) in Dalarna got a higher percentage of the absentee ballots. Since (MP) got one more permanent seat in step (ii), they were given one adjustment seat less in step (iv). If the total proportionality of step (iii) is not changed this would normally, i.e. if the BUT had not happened, give an extra adjustment to the party that lost a permanent seat in this case (S). But since (S) already had too many seats the adjustment seat instead went to a third party Centerpartiet (C). They got this new seat in the constituency of S\"odermanland. Thus the higher percentage of votes for (MP) caused the transfer of a seat from a coalition partner (S) to a non-coalition partner (C). This was most likely not according to the voters' intentions and clearly an unwanted behavior of an electoral system.

The same principles are used also in the election to the 20 county councils in Sweden. 
In 2010 the BUT happened in
9 out of 20 county councils. One reason it was more likely to happen in 2010 than previously is that the number 
of parties above the 4\%-threshold increased to 8. Similarly for the county councils, where there is a 3\%-
threshold, the number of parties has increased to between 6 and 9.

%%%%%%%%%%%%%%%%
\section{Possible cures and dynamic adjustment} \label{S:Dynamic}
The deviation from proportionality might not be considered very large, but could be of large political 
significance. It is also possible to improve the system to
better fulfill the intention of the law makers with rather small technical adjustments.
There are at least three obvious ways of lowering the probability of the BUT to happen:

\smallskip\noindent
1. Fewer constituencies\footnote{Possibility 1 also has the drawback that it will have to be repeated separately for each county council, whereas possibilities 2,3 and 4 easily carry over directly.}.  (probably politically very difficult to achieve)\\
2. Changing the divisor $1.4$. (around $1.2$ seems to be better by some simulations\footnote{Simulations mostly on national level to Riksdagen.}, \cite{Lanke, Fr}) We see no theoretical reason to think that some particular divisor would give the minimum number of adjustment seats, but heuristically a lower divisor would make it easier for small parties to get their first seat in a constituency and thereby lowering the probability that the big parties get too many seats in total. But too low a divisor would on the other hand give increased risk of a small party getting too many permanent seats.\\
3. More adjustment seats and fewer permanent seats.\\

 The third possibility is perhaps the most attractive of the three above, but it has the drawback that to be reasonably  certain
 to have proportionality between parties, one would need to use a rather large number of adjustment seats. 
 This could cause a not so good distribution of seats between constituencies. In particular one does not want the 
 small constituencies to be unfavourably treated. We would like here to draw attention to a fourth possibility, 
 which would mean a somewhat larger change to the system by using what we call the 
{\bf dynamic} method. \\

4. (Dynamic adjustment) Let the adjustment seats be exactly as many as needed to give national proportionality between parties in that particular election\footnote{Because of the Alabama paradox for Hamilton's method, it would be preferable to change to \web{} also for distribution of permanent seats between constituencies, see discussion in \refS{S:Implement}.}. That is, we keep the total number of seats fixed but the point at which one switches from permanent seats to adjustment seats is flexible.

\medskip
The advantages of dynamic adjustment are several:

\begin{romenumerate}
\item{} It always gives proportionality between parties.
\item{} It is a rather small change to the current system.
\item{} It does not use more adjustment seats than necessary and thus increases the likelihood of good geographic representation, see \refS{S:Simulate}.
\item{} It is easy to implement, see \refS{S:Implement}. 
\item{} It avoids the strange effects in the current system as described in \refS{S:Sweden}.
\end{romenumerate}

The most obvious advantage with the dynamic method compared to the current method is that it will give
 proportionality between parties. It is our guess that under most circumstances the number of 
adjustment seats will be relatively low, see \refS{S:Simulate}, which would imply that no constituency would 
get unreasonably few seats.

Another advantage is that we avoid the troubles that the current Swedish system runs into when the number of
adjustment seats needs to be more than 39, i.e. when the BUT occurs. For instance, as explained at the end
of \refS{S:Sweden}, extra votes on party A might move a seat from party B to party C in the current system.
This resembles the so called population paradox of Hamilton's method (a.k.a. 
the method of largest remainders), see e.g. \cite{BY}. One reason to use the \web{} in the first place is to 
avoid such things to happen. This cannot happen with the dynamic method as the total number of seats for any 
party is given by a calculation on the national level.

We have found that the negative sides of dynamic adjustment are, compared to the current Swedish system, as follows.

\begin{itemize}
\item[(vi)] It is not possible to give a theoretical upper bound to the number of adjustment seats needed in the worst case. 
\item[(vii)] It could be non-monotone for a single candidate, see \refS{S:Non-monotone}.
\item[(viii)]  The details are slightly more complicated to explain to the general public.
\end{itemize}

\smallskip
The lack of a theoretical upper bound on the number of adjustment seats could, at least in a constructed example, lead to the smallest constituency (Gotland) getting only one or perhaps even zero seats. This must clearly be avoided. Such a scenario is unlikely and in a practical implementation one can 
easily add an extra requirement that the number of adjustment seats is at most some fixed number, see 
\refS{S:Implement}. The non-monotonicity for single candidates is more cumbersome, see discussion in 
\refS{S:Non-monotone}.

A proposed system with a reminiscent dynamic adjustment was discussed in \cite{BF}, but it is rather different. For instance,
there the number of permanent seats is kept fixed (one in each constituency) 
and then the total size of the parliament is allowed to vary.

\section{Aspects on an implementation} \label{S:Implement}
No matter which formulation of the \web{} is used it is easy to give a direct description of the method of dynamic adjustment.
Let us describe in some detail how this could be done based on the formulation used in \refS{S:Sweden}.
With this formulation the seats are distributed one at a time to the entity most 
deserving the next seat.
Thus the outcome is not only a total of how many seats each party (constituency) should have but also a list 
$L$ in 
which order they receive the seats.
An outline of the dynamic method with some comments is as follows. Let $m$ be the total number of seats in 
the parliament. Note that the total number of seats is fixed.

\begin{enumerate}
\item{} Before election the $m$ seats are given a tentative distribution to the constituencies with the
\web{}  (the odd number method). 
This gives an ordered list, called $L$ with the order in which the seats are to be given to the constituencies. 
\item{} Voting takes place.
\item{} A proportional distribution (the \web) of the seats between the parties is created. Let $r_j$ be 
the number of seats party $j$ should have, i.e. the column sum of column $j$ in the output matrix $M$.\footnote{In this step 
we assume that the \web{} is used again, but other methods are also possible. Also 
other constraints such as a minimum threshold for small parties may be used in this step.}
\item{} \label{stop} We now follow the order of the constituencies given by $L$. 
If the next seat should be given to constituency
$i$ then an extra seat is given to row $j$ of the matrix $M$. Within the constituency (each row) 
the seat is distributed  with the \web{} between parties. This process continues until
some party $j$ is just about to get a total of seats larger than $r_j$, that is the column sum of column $j$ is about to 
be exceeded. Then this part of the process ends before that seat is given out.\footnote{To implement the 12\%-
threshold, see \refS{S:Sweden}, one has to take some care. One 
could let a party $A$ with more than 12\% of the votes in a
constituency but fewer than 4\% nationwide take part in the distribution of permanent seats in this step, but not 
in the distribution of adjustment seats in the next step. If party $A$ gets a seat the total number of seats has to 
be reduced for the party getting the last seat in step (3), say party $B$. If that party already has been given all 
its seats then party $B$ gets to keep its seats and the process stops just before party $A$ is given any seat.} We 
say that the {\bf stop for permanent seats} has occurred.
\item{} The remaining seats are distributed as adjustment seats as in \refS{S:Sweden}. That is to say each column 
in $M$ will get the correct column sum and the adjustment seats will be distributed within the column 
using the \web{} with initial values given by the seats already handed out in the previous step.
\end{enumerate}

\begin{remark}
With this method, we guarantee proportionality between parties, since we have made the column sums  
unnegotiable side constraints. The number of adjustment seats will be  exactly as small as possible,
but still giving the wanted total outcome. 
\end{remark}
\begin{remark}
It is not a key point point in the method how the list $L$ is generated in (1). 
The nice feature of the \web{} here is that the length of the list does not influence the ordering in the list, which implies that
the implementation described above gives the same result as minimizing the number of adjustment seats. This is not true for 
Hamilton's method  because of the Alabama paradox, see e.g. \cite{BY, JL}.
\end{remark}

If we are afraid of the number of adjustment seats being too large, the law makers could decide on an a
minimum $m'$ of permanent seats. Then step \eqref {stop} has to be changed to 

\begin{itemize}
\item[(\ref{stop}')] We now follow the order of the constituencies given by $L$. 
If the next seat should be given to constituency
$i$ then an extra seat is given to row $i$ of the matrix $M$. Within each constituency (each row) 
the seat is distributed by the \web{} between parties. This process continues until
some party $j$ is just about to get a total of seats larger than $r_j$, that is, the column sum of column $j$ is about to 
be exceeded. {\em If already at least $m'$ permanent seats have been distributed}, this part of the process ends 
before the seat causing the exceedance is distributed. {\em Otherwise, this step continues until exactly $m'$ permanent seats have been 
distributed.}
\end{itemize}

Yet another possibility would be to start the process by giving every constituency 2 permanent seats, 
to prevent the, unlikely, event that the smallest constituency Gotland only got one seat.

The reader might ask why after the stop for permanent seats occurs, one does not simply remove that party and iterate (4) with the remaining parties. The reason is that the very last seat will then be bound to go to a certain party in a certain constituency, where it might have very few votes or even no votes at all.

\section{Theoretical weaknesses of dynamic adjustment}\label{S:theory}
In this section we discuss three possible drawbacks of the dynamic method from a theoretical viewpoint.

\subsection{Number of adjustment seats} \label{S:Adjust}

One drawback with the dynamic method is that one cannot really give any good theoretical bounds on how 
many adjustment seats that might be needed. As the following example suggests the adjustment 
seats could in theoretically constructed examples be very many, in fact almost all seats.

\smallskip
\noindent
{\bf Example 1}\\
Assume we have two parties A and B. The total number of seats is 208. Party A gets 63000 votes and party B gets 3010. This implies that A should have 199 and B 9 seats.
There are 110 constituencies and the votes are distributed as follows. We assume 100\% voter turnout in these theoretical examples.

\medskip
\begin{tabular}{| l | r| r |}
\hline
Party & A & B \\
\hline
Constituency 1-10 & 300 & 301 \\
Constituency 11-110 & 600 & 0 \\
\hline
Total & 63 000 & 3 010\\
\hline
\end{tabular}

\smallskip
Clearly Constituencies 1-10 will be first in our list $L$ when we distribute seats. Party
B will get a seat in each of these and hence the stop for permanent seats will occur just before giving B a tenth seat. Then the remaining 199 seats will be given to A as adjustment seats. This could easily be generalized by
increasing the number of constituencies of the second kind.
\smallskip

However, note in this example that the overall distribution of seats is very reasonable. Thus the fact that the 
number of adjustment seats is large does not imply that we will have a poor geographic proportionality. Also 
note that the main characteristic is a large number of small constituencies where one party gets all the votes. 
This is far from the reality in Sweden. In Table \ref{T:adjust_seats} we present the number of adjustment seats that would have been 
needed in the elections 1970-2010. In Figure \ref{F:adj} is a diagram of the number of adjustment seats needed in
 10000 simulations done around the election result of 2010.

\medskip

\subsection{Monotonicity for parties but not for candidates} \label{S:Non-monotone}
A worse problem is the theoretical possibility that a vote for a certain party might not be to the advantage of 
a candidate for the party running in the constituency where the vote was cast. 

\smallskip
\noindent
{\bf Example 2}\\
Consider the following situation. Two parties are competing for three seats in a total of three constituencies.
Assume first that the outcome of the election is as follows.

\medskip
\begin{tabular}{|l | r| r || r |}
\hline
Constituency$\backslash$  Party &A & B& Total\\
\hline
I & 97 & 98 &  195 \\
II & 101 & 100  & 201 \\
III & 102 & 101 & 203 \\
\hline
Total& 300 & 299  & 599\\
\hline
\end{tabular}

\smallskip
Party A should by the \web{} get two of the three seats while Party B should get one. The list $L$
gives one seat to each constituency in the order III, II and  I. Thus Party A will win the two seats in constituencies 
III and II, Party B the seat in I and no adjustment seats would be necessary.

But if one voter in constituency I had changed his mind in the last minute and decided to vote for Party B instead of Party A we would get the following table of votes.

\medskip
\begin{tabular}{|l | r| r || r |}
\hline
Constituency$\backslash$  Party &A & B& Total\\
\hline
I & 96 & 99 &  195 \\
II & 101 & 100  & 201 \\
III & 102 & 101 & 203 \\
\hline
Total& 299 & 300  & 599\\
\hline
\end{tabular}

\smallskip
With this outcome Party B should get two seats in total and Party A only one. The list $L$ is not changed so 
Party A gets first the seat in III and then the stop occurs and Party B gets both its seats as adjustment seats, 
which will be in constituencies III and II. Thus the candidate for Party B in constituency I lost her seat because of
the extra vote she got. 
The candidate lost it to two other members of Party B, and in fact, in the dynamic 
method it cannot be bad
for a party to gain extra votes. But for a single candidate it could be as Example 2 shows.  In Section \ref{S:SimNon-monotone}
we test how likely this is to happen, by simulations close to the election result from 2010.
This effect is not possible in the current 
Swedish system.

\subsection{Non-monotonicity for constituencies} \label{S:Const}
A problem with the current Swedish system that, as far as we know, has not been discussed at all is the fact 
that a high voter turnout in a constituency might lead to fewer seats for that constituency.
Let us for instance look at the constituency of Halland in the 2006 national election in Sweden. 
The numbers of votes for the seven parties represented are given in Table \ref{F:Halland}. Halland was 
unlucky and got none of the 39 adjustment seats in 2006. If (KD) had 
received 1337 fewer votes, they would not have received a permanent seat in Halland. A look at the national 
number of votes shows that the total number of seats to each party would not have changed and (KD)
would instead have gotten an adjustment seat in Halland.

\begin{table}[htb]
\begin{tabular}{|l | r| r | r || r | r | r | r |}
\hline
Party &	actual & perm.  & adj. & assumed &
new & new perm. & new adj. \\ 
 &	votes & seats & seats & change &
votes& seats & seats\\ 
\hline
M	&53257	&3	&0 &0 &53257 &4	&0\\
\hline
C	&18589	&1	&0 &0 &18589 &1	&0\\
\hline
FP	&13798	&1	&0 &0 & 13798 &1	&0\\
\hline
KD	&11987	&1	&0 &-1337 & 10650  &0	&1\\
\hline
S	&56747	&4	&0 &0 &56747&4	&0\\
\hline
VP	&7110	&0	&0 &0 &7110&0	&0\\
\hline
MP	&7236	&0	&0 &0 &7236&0	&0\\
\hline	
\end{tabular}
\caption{If some voters in Halland had not voted in 2006 election, the constituency of Halland could have  received one more seat in the Swedish parliament. With the 1337 {\it fewer} votes Halland would have received 11 seats instead of 10.} \label{F:Halland}
\end{table}

It is easy to find such examples, we could have chosen several other constituencies in the same election. 
In particular the fact that the quotient $1.4$ is used when determining 
the permanent seats in each region but not for adjustment seats makes this rather likely to happen. 

This effect relies on the fact that there are two steps in the assignment method, and would be possible also with the dynamic adjustment method.

\medskip
One subtle point of possible interest is that in the current Swedish system the number of permanent seats to a
constituency is determined by the number of entitled voters, whereas the adjustment seats are distributed
based on the number of votes given. This gives a certain mix of relevance to these two variables.
In the dynamic system this mix would not be determined beforehand, but rather as a result of
the algorithm, which might be considered unsatisfactory from a theoretical viewpoint. The voter turnout did not vary much in last election, from $79.63\%$ to $86.36\%$.

\section{Methods for measuring disproportionality} \label{errorMeasures}
\subsection{Disproportionality for parties, constituencies and voter representation}
When comparing the vote result with a seat distribution, there are three types of disproportionality to consider: the most basic one refers to voter representation, i.e., how well the seat distribution represents the voters. If 1\% of the votes is for party A in constituency I, it would be good if (close to) 1\% of the seats were assigned to (party A, constituency I).
In Swedish tradition, two other views have been more in focus: party error and constituency error. The constituency error is the disproportionality when comparing the number of seats assigned to each constituency compared with the number of votes in that constituency. An alternative definition here would be to compare the number of seats for each constituency with the number of entitled voters in that constituency. In this paper we will however measure the disproportionality with respect to the actual vote result. The party error is the disproportionality when comparing  the number of seats assigned to each party compared with the number of votes for that party.

\smallskip
Note that much of the reported problems with disproportionality in the 2010 elections in Sweden is due to the way that constituency error is being bounded.
Among all proposed methods, it seems that the mathematically more advanced method 
called the biproportional method, see e.g \cite{BP}, is best at
lowering both the party error and the constituency error simultaneously. The biproportional method could be argued for as the best method to solve the problems in Sweden. The main reasons legislatures would prefer the dynamic method is, as far as we can see, that they want to do only a small change to the current system and secondly that the more advanced mathematics needed makes it less accessible to the general public in full detail.

\subsection{Quantifying disproportionality}
Several different methods for measuring disproportionality have been suggested. We present results here in two different ways. The first way is by a normalized $l^1$-method suggested by Loosemore and Hanby \cite{LH}: the measure $E_{LH}$ of disproportionality for the seat representation $s=\{s_i\}$ relative to the vote result $v=\{v_i\}$ is
\begin{equation}
E_{LH} = 50 \sum |(v_i/V) - (s_i/S)|,
\end{equation}
where the summation is performed over all parties or all constituencies or all pairs of (party, constituency) 
depending on the error category, and where $S$ is the sum of seats and $V$ the sum of votes for parties entering the parliament.
The number 50 comes  from the fact that they used percentages and divided the results with 2.

Most commonly known divisor methods minimize some measure of deviation from exact proportionality. 
One can argue, see Gallagher \cite{G91}, that thus each such method comes with its own taylored measure of disproportionality. Since the current
Swedish system is based on the \web, we have chosen as a second measure to use the method connected to the \web{}, which we here will call the SL-measure. It is relative to the size of the party (constituency) so that if a party 
with $5\%$ of the votes receives $6\%$ of the seats is considered as larger deviation than if a party 
with $10\%$ of the votes receives $11\%$ of the seats.  
The formula for the SL-measure is, with notation as above 
\begin{equation}
E_{SL} =100 \sum \frac{1}{v_i/V}\left((v_i/V) - (s_i/S)\right)^2,
\end{equation}
where again the summation is performed over all parties or all constituencies.

\section{Back-testing} \label{S:backtesting} In this section we analyze how the dynamic method would have performed in some recent national elections in Sweden. In particular with reference to the theoretical concerns discussed in Section \ref{S:theory} on how many adjustment seats would have been needed, disproportionality and the possibility of non-monoticity for candidates.
 
\subsection{Data}
The vote count from the national parliament  election in Sweden for the elections $\{1973, 1976, 1979, 1982, 1985, 1988, 1994, 1998, 2002, 2006, 2010\}$ have been obtained from Valmyndigheten and SCB\footnote{We thank Jan Lanke that has been kind enough to supply us with the values for 1970 and 1991 in Tables \ref{T:adjust_seats} and \ref{T:LHparty}.}. 
We were not able to retrieve the complete election data for the election of 1991.

\subsection{Number of adjustment seats}
In Table \ref{T:adjust_seats} we give how many adjustment seats the dynamic method would have used for the outcome
of the Swedish national parliament elections. The second row gives the number of adjustment seats needed 
for the modified \web{} to give proportionality. The total number of seats was changed in the early seventies.

\begin{table}[h] %\label{T:BackAdj}
\begin{tabular}{|l | r | r | r | r | r | r | r| r | r | r | r | r | r |}
\hline
Year & \hspace{-3pt}1970\hspace{-3pt} &\hspace{-3pt}1973\hspace{-3pt} &\hspace{-3pt}1976\hspace{-3pt} &\hspace{-3pt}1979\hspace{-3pt} &\hspace{-3pt}1982\hspace{-3pt} &\hspace{-3pt}1985\hspace{-3pt} &\hspace{-3pt}1988\hspace{-3pt} &\hspace{-3pt}1991\hspace{-3pt}  &\hspace{-3pt}1994\hspace{-3pt} &\hspace{-3pt}1998\hspace{-3pt} &\hspace{-3pt}2002\hspace{-3pt} &\hspace{-3pt}2006\hspace{-3pt} &\hspace{-3pt}2010\hspace{-3pt} \\
\hline
Dynamic, $1.0$& 9\hspace{-3pt} & 18& 3& 17& 38& 1& 21& 10& 48& 37& 32& 42& 52\\
Dynamic, $1.4$& 19\hspace{-3pt} & 17& 16& 26& 25& 12& 51& 10& 36& 36& 30& 30& 57\\
\hline
Parties& 5\hspace{-3pt} & 5& 5& 5& 5& 5& 6& 7 & 7& 7& 7& 7& 8\\
\hline
Total $\#$ of seats \hspace{-6pt} & 350\hspace{-3pt} & 350 & 349 & 349& 349& 349& 349& 349 & 349& 349& 349& 349& 349\\
\hline
\end{tabular}
\vskip2mm
\caption{Number of adjustment seats needed to give proportionality for the dynamic method, with the pure and 
modified \web{} respectively.}\label{T:adjust_seats}
\end{table}
%\medskip
\subsection{Disproportionality comparison between existing and dynamic adjustment method}
The party disproportionality between the current Swedish system and the dynamic system will be the same every time the number of adjustment seats has 
been enough, i.e. the BUT has not occurred. 

{}From Table \ref{T:LHparty} below we see that the party disproportionality is in general small. 
But, as noted in BUT, for the elections in 2010 and 1988, they are not merely rounding errors for the current method. 
In agreement with the design, the dynamic method 
only has a party error originating from rounding. 

\begin{figure}[ht]
\begin{tikzpicture}[xscale=0.4,yscale=0.77]
\draw [<->] (70,4.4) -- (70,0) -- (104,0);
%\draw[help lines] (1,0) grid (4,2);
\node[red] at (73,1.81) {$\otimes$};
\node[red]  at (76,1.23) {$\otimes$};
\node[red]  at (79,1.82) {$\otimes$};
\node[red]  at (82,2.76) {$\otimes$};
\node[red]  at (85,1.18) {$\otimes$};
\node[red]  at (88,1.81) {$\otimes$};
\node[red]  at (91,3.26) {$\otimes$};
\node[red]  at (94,3.01) {$\otimes$};
\node[red]  at (97,2.13) {$\otimes$};
\node[red]  at (100,3.32) {$\otimes$};
\node[red]  at (103,3.49) {$\otimes$};

\node[below] at (73,1.81) {\small$1.81$};
\node[below]  at (76,1.23) {\small$1.23$};
\node[below]  at (79,1.82) {\small$1.82$};
\node[above]  at (82,2.76) {\small$2.76$};
\node[below]  at (85,1.18) {\small$1.18$};
\node[below]  at (88,1.81) {\small$1.81$};
\node[below]  at (91,3.26) {\small$3.26$};
\node[below]  at (94,3.01) {\small$3.01$};
\node[above]  at (97,2.13) {\small$2.13$};
\node[below]  at (100,3.32) {\small$3.32$};
\node[below]  at (103,3.49) {\small$3.49$};

\node [ blue] at (73,2.34) {$\bigstar$};
\node [ blue] at (76,2.65) {$\bigstar $};
\node [ blue]at (79,2.28) {$\bigstar $};
\node [ blue]at (82,2.68) {$\bigstar $};
\node [ blue] at (85,1.65) {$\bigstar $};
\node [ blue] at (88,2.73) {$\bigstar $};
\node [ blue] at (91,3.79) {$\bigstar $};
\node [ blue] at (94,3.02) {$\bigstar $};
\node [ blue] at (97,2.00) {$\bigstar $};
\node [ blue] at (100,3.38) {$\bigstar $};
\node [ blue] at (103,3.75) {$\bigstar $};

\node [ above] at (73,2.34) {\small$2.34$};
\node [ above] at (76,2.65) {\small$2.65$};
\node [ above]at (79,2.28) {\small$2.28$};
\node [ below]at (82,2.68) {\small$2.68$};
\node [ left] at (85,1.65) {\small$1.65$};
\node [ above] at (88,2.73) {\small$2.73$};
\node [ above] at (91,3.79) {\small$3.79$};
\node [ above] at (94,3.02) {\small$3.02$};
\node [ left] at (97,2.00) {\small$2.00$};
\node [ above] at (100,3.38) {\small$3.38$};
\node [ above] at (103,3.75) {\small$3.75$};

\node [ below] at (73,-0.3) {$1973$};
\node [ below] at (76,-0.3) {1976};
\node [ below] at (79,-0.3) {1979};
\node [ below] at (82,-0.3) {1982};
\node [ below] at (85,-0.3) {1985};
\node [ below]at (88,-0.3) {$1988$};
\node [ below]at (91,-0.3) {$1994$};
\node [ below] at (94,-0.3) {$1998$};
\node [ below] at (97,-0.3) {$2002$};
\node [ below] at (100,-0.3) {$2006$};
\node [ below] at (103,-0.3) {$2010$};
\node[right] at (70,4){LH};
\end{tikzpicture}

\caption{LH measure of constituency disproportionality. Blue stars $\bigstar$ are the current Swedish system and red circles with a cross $\otimes$ are what the dynamic method would have given.}\label{F:LH_Const}
\end{figure}
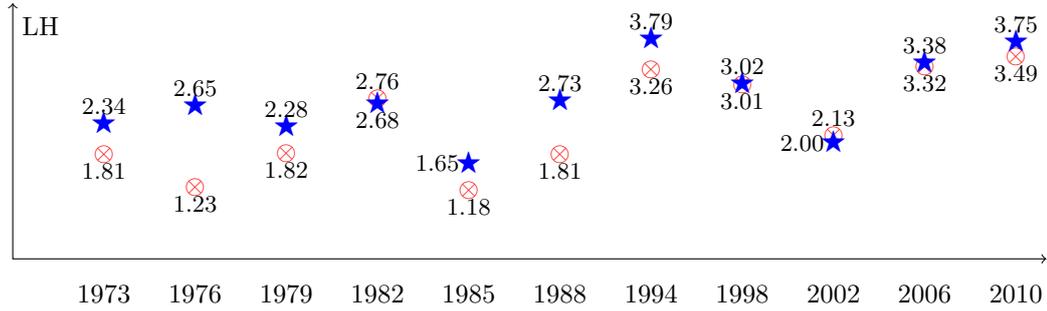

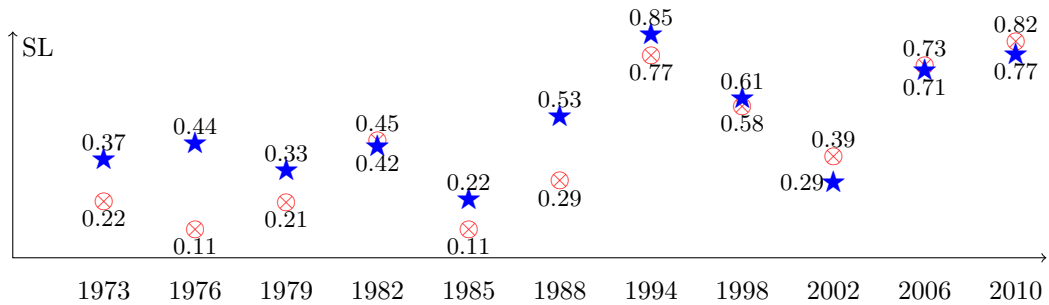
\begin{figure}[htb]
\begin{tikzpicture}[xscale=0.4,yscale=1]
\draw [<->] (70,3) -- (70,0) -- (104,0);
%\draw[help lines] (1,0) grid (4,2);
\node [red] at (73,0.75) {$\otimes$};
\node [red] at (76,0.38) {$\otimes$};
\node  [red]at (79,0.74) {$\otimes$};
\node  [red]at (82,1.56) {$\otimes$};
\node [red] at (85,0.38) {$\otimes$};
\node [red] at (88,1.02) {$\otimes$};
\node [red] at (91,2.68) {$\otimes$};
\node [red] at (94,2.01) {$\otimes$};
\node [red] at (97,1.35) {$\otimes$};
\node [red] at (100,2.55) {$\otimes$};
\node [red] at (103,2.87) {$\otimes$};

\node[below] at (73,0.75) {\small$0.22$};
\node[below]  at (76,0.38) {\small$0.11$};
\node[below]  at (79,0.74) {\small$0.21$};
\node[above]  at (82,1.56) {\small$0.45$};
\node[below]  at (85,0.38) {\small$0.11$};
\node[below]  at (88,1.02) {\small$0.29$};
\node[below]  at (91,2.68) {\small$0.77$};
\node[below]  at (94,2.01) {\small$0.58$};
\node[above]  at (97,1.35) {\small$0.39$};
\node[above]  at (100,2.55) {\small$0.73$};
\node[above]  at (103,2.87) {\small$0.82$};

\node [ blue] at (73,1.31) {$\bigstar$};
\node [ blue] at (76,1.52) {$\bigstar$};
\node [ blue]at (79,1.16) {$\bigstar$};
\node [ blue]at (82,1.48) {$\bigstar$};
\node [ blue] at (85,0.77) {$\bigstar$};
\node [ blue] at (88,1.87) {$\bigstar$};
\node [ blue] at (91,2.96) {$\bigstar$};
\node [ blue] at (94,2.11) {$\bigstar$};
\node [ blue] at (97,1.00) {$\bigstar$};
\node [ blue] at (100,2.48) {$\bigstar$};
\node [ blue] at (103,2.70) {$\bigstar$};

\node [ above] at (73,1.31) {\small$0.37$};
\node [ above] at (76,1.52) {\small$0.44$};
\node [ above]at (79,1.16) {\small$0.33$};
\node [ below]at (82,1.48) {\small$0.42$};
\node [ above] at (85,0.77) {\small$0.22$};
\node [ above] at (88,1.87) {\small$0.53$};
\node [ above] at (91,2.96) {\small$0.85$};
\node [ above] at (94,2.11) {\small$0.61$};
\node [ left] at (97,1.00) {\small$0.29$};
\node [ below] at (100,2.48) {\small$0.71$};
\node [ below] at (103,2.70) {\small$0.77$};

\node [ below] at (73,-0.2) {$1973$};
\node [ below] at (76,-0.2) {1976};
\node [ below] at (79,-0.2) {1979};
\node [ below] at (82,-0.2) {1982};
\node [ below] at (85,-0.2) {1985};
\node [ below]at (88,-0.2) {$1988$};
\node [ below]at (91,-0.2) {$1994$};
\node [ below] at (94,-0.2) {$1998$};
\node [ below] at (97,-0.2) {$2002$};
\node [ below] at (100,-0.2) {$2006$};
\node [ below] at (103,-0.2) {$2010$};
\node[right] at (70,2.8){SL};
\end{tikzpicture}

\caption{SL measure of constituency disproportionality. Blue stars $\bigstar$ are the current Swedish system and red circles with a cross $\otimes$ are what the 
dynamic method would have given.}\label{F:SL_Const}
\end{figure}

\begin{table}[h]
\begin{tabular}{|l | r | r | r | r | r | r | r| r | r | r | r | r | r |}
\hline
Year & \hspace{-3pt}1970\hspace{-3pt} &\hspace{-3pt}1973\hspace{-3pt} & \hspace{-3pt}1976\hspace{-3pt} & \hspace{-3pt}1979\hspace{-3pt} & \hspace{-3pt}1982\hspace{-3pt} & \hspace{-3pt}1985\hspace{-3pt} & \hspace{-3pt}1988\hspace{-3pt} & \hspace{-3pt}1991\hspace{-3pt}  & \hspace{-3pt}1994\hspace{-3pt} & \hspace{-3pt}1998\hspace{-3pt} & \hspace{-3pt}2002\hspace{-3pt} & \hspace{-3pt}2006\hspace{-3pt} & \hspace{-3pt}2010\hspace{-3pt} \\
\hline
Current method&	0.19\hspace{-3pt} & 0.10\hspace{-3pt} &	0.20\hspace{-3pt} &	0.20\hspace{-3pt} &	0.21\hspace{-3pt} &	0.11\hspace{-3pt} &	0.31\hspace{-3pt} & 0.24\hspace{-3pt} &	0.27\hspace{-3pt} &	0.20\hspace{-3pt} &	0.22\hspace{-3pt} &	0.29\hspace{-3pt} &	1.15\hspace{-3pt} \\
\hline
Dynamic method &	0.19\hspace{-3pt} & 0.10\hspace{-3pt} &	0.20\hspace{-3pt} &	0.20\hspace{-3pt} &	0.21\hspace{-3pt} &	0.11\hspace{-3pt} &	0.18\hspace{-3pt} & 0.24\hspace{-3pt} &	0.27\hspace{-3pt} &	0.20\hspace{-3pt} &	0.22\hspace{-3pt} &	0.29\hspace{-3pt} &	0.24\hspace{-3pt} \\
\hline
\end{tabular}
\vskip2mm
\caption{Measure LH of party disproportionality.}\label{T:LHparty}
\end{table}
Less easy to predict is the comparison between 
the two methods for constituency disproportionality. In Figure \ref{F:LH_Const} and
Figure \ref{F:SL_Const} we give the disproportionality for the current method compared to the dynamic with the LH and the SL measure respectively. The disproportionality is not very large and as can be seen from the diagrams, the dynamic method has a slightly smaller measure of disproportionality. We think that this is an effect of the current Swedish method using 
the modified \web, i.e. the $1.4$-factor, which gives an unwanted disproportionality, see 
discussion in \refS{S:Simulate}.

\section{Simulations} \label{S:Simulate}
\subsection{Perturbations around 2010 election results}
In order to investigate both the behaviour of the proposed method and the current method for results similar to the latest election result, we calculate the outcome for these methods for simulated election results close the outcome of the latest election.

The simulation is done in the following way. For party $j$ in constituency $i$, the simulated result $v^s_{ij}$ is set to 
\begin{equation}\label{perturbedSim}
v^s_{ij} = round(v_{ij} * p_j*x_{ij}),
\end{equation}
where the $p_j,x_{ij}$ are independent stochastic variables with uniform distribution in $(0.9, 1.1)$ and $v_{ij}$ is the actual number of votes in the 2010 election. Here $p_j$ is thought of as a 
nationwide change for each party and $x_{ij}$ as local fluctuations. 
We create $N=10000$ simulated election results according to \eqref{perturbedSim}, and record the 
number of adjustment seats, the various disproportionality measures and check if there is a threat of non-monotonicity for 
a candidate, as described in Section \ref{S:Non-monotone}.

\subsection{Number of adjustment seats}

\begin{figure}[htb]
\begin{tikzpicture}[xscale=0.13,yscale=0.055]
\draw [<->] (0,50) -- (0,0) -- (100,0);
\draw [fill=gray] (20,0) rectangle (29,0);
\node [below] at (25,0) {20-29};
\node [above ] at (25,0) {0};
\draw [fill=gray] (30,0) rectangle (39,4.64);
\node [below] at (35,0) {30-39};
\node [above] at (35,4.64) {464};
\draw [fill=gray] (40,0) rectangle (49,48.92);
\node [below] at (45,0) {40-49};
\node [above] at (45,48.92) {4892};
\draw [fill=gray] (50,0) rectangle (59,37.56);
\node [below] at (55,0) {50-59};
\node [above] at (55,37.56) {3756};
\draw [fill=gray] (60,0) rectangle (69,8.64);
\node [below] at (65,0) {60-69};
\node [above] at (65,8.64) {864};
\draw [fill=gray] (70,0) rectangle (79,0.24);
\node [below] at (75,0) {70-79};
\node [above] at (75,0.24) {24};
\draw [fill=gray] (80,0) rectangle (89,0);
\node [below] at (85,0) {80-89};
\node [above] at (85,0) {0};
\end{tikzpicture}
\caption{Histogram over the number of adjustment seats needed for proportionality in 10000 simulations of the 
dynamic method with the modified \web. Average $49.6$, maximum $74$.}\label{F:adjmod}
\end{figure}
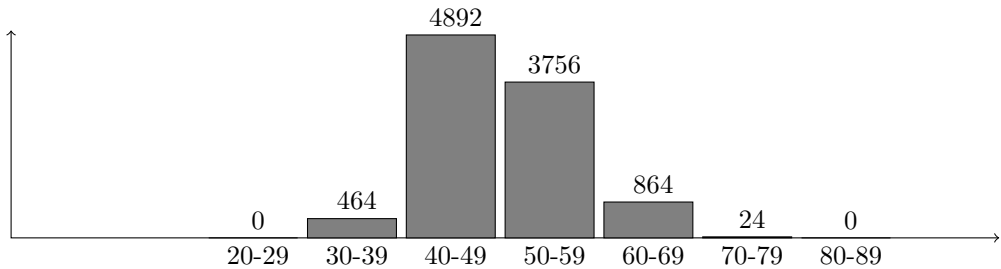

The numbers of adjustment seats needed in the simulations are depicted in Figures \ref{F:adj} and \ref{F:adjmod}. The average 
when using the standard \web{} was $52.3$ and largest 107, when using the modified \web{} the average was $49.6$  and the 
largest 74. The standard deviation seems to be smaller when using the modified version.

\begin{figure}[hb] 
\begin{tikzpicture}[xscale=0.13,yscale=0.09]
\draw [<->] (0,30) -- (0,0) -- (110,0);
\draw [fill=gray] (20,0) rectangle (29,4.65);
\node [below] at (25,0) {20-29};
\node [above ] at (25,4.65) {465};
\draw [fill=gray] (30,0) rectangle (39,15.68);
\node [below] at (35,0) {30-39};
\node [above] at (35,15.68) {1568};
\draw [fill=gray] (40,0) rectangle (49,20.82);
\node [below] at (45,0) {40-49};
\node [above] at (45,20.82) {2082};
\draw [fill=gray] (50,0) rectangle (59,28.08);
\node [below] at (55,0) {50-59};
\node [above] at (55,28.08) {2808};
\draw [fill=gray] (60,0) rectangle (69,19.52);
\node [below] at (65,0) {60-69};
\node [above] at (65,19.52) {1952};
\draw [fill=gray] (70,0) rectangle (79,8.41);
\node [below] at (75,0) {70-79};
\node [above] at (75,8.41) {841};
\draw [fill=gray] (80,0) rectangle (89,2.34);
\node [below] at (85,0) {80-89};
\node [above] at (85,2.34) {234};
\draw [fill=gray] (90,0) rectangle (99,0.43);
\node [below] at (95,0) {90-99};
\node [above] at (95,0.43) {43};
\draw [fill=gray] (100,0) rectangle (109,0.07);
\node [below] at (105,0) {100-109};
\node [above] at (105,0.07) {7};
\end{tikzpicture}
\caption{Histogram over the number of adjustment seats needed for proportionality in 10000 simulations of 
the dynamic method with the standard \web. Average $52.3$, maximum $107$.}\label{F:adj}
\end{figure}
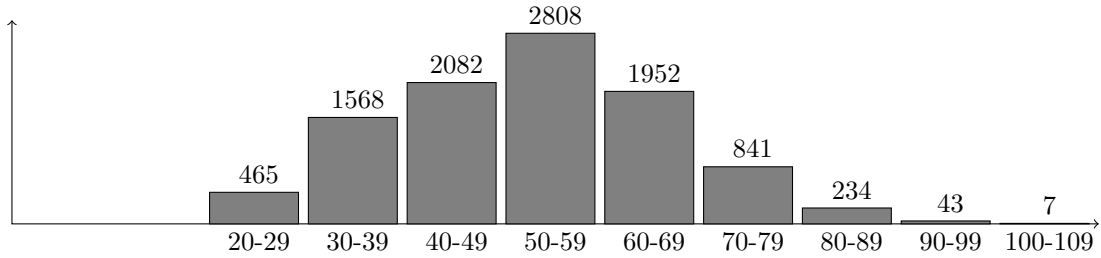

\subsection{Proportionality for constituencies}

\begin{figure}[htb]
\includegraphics[height=6cm]{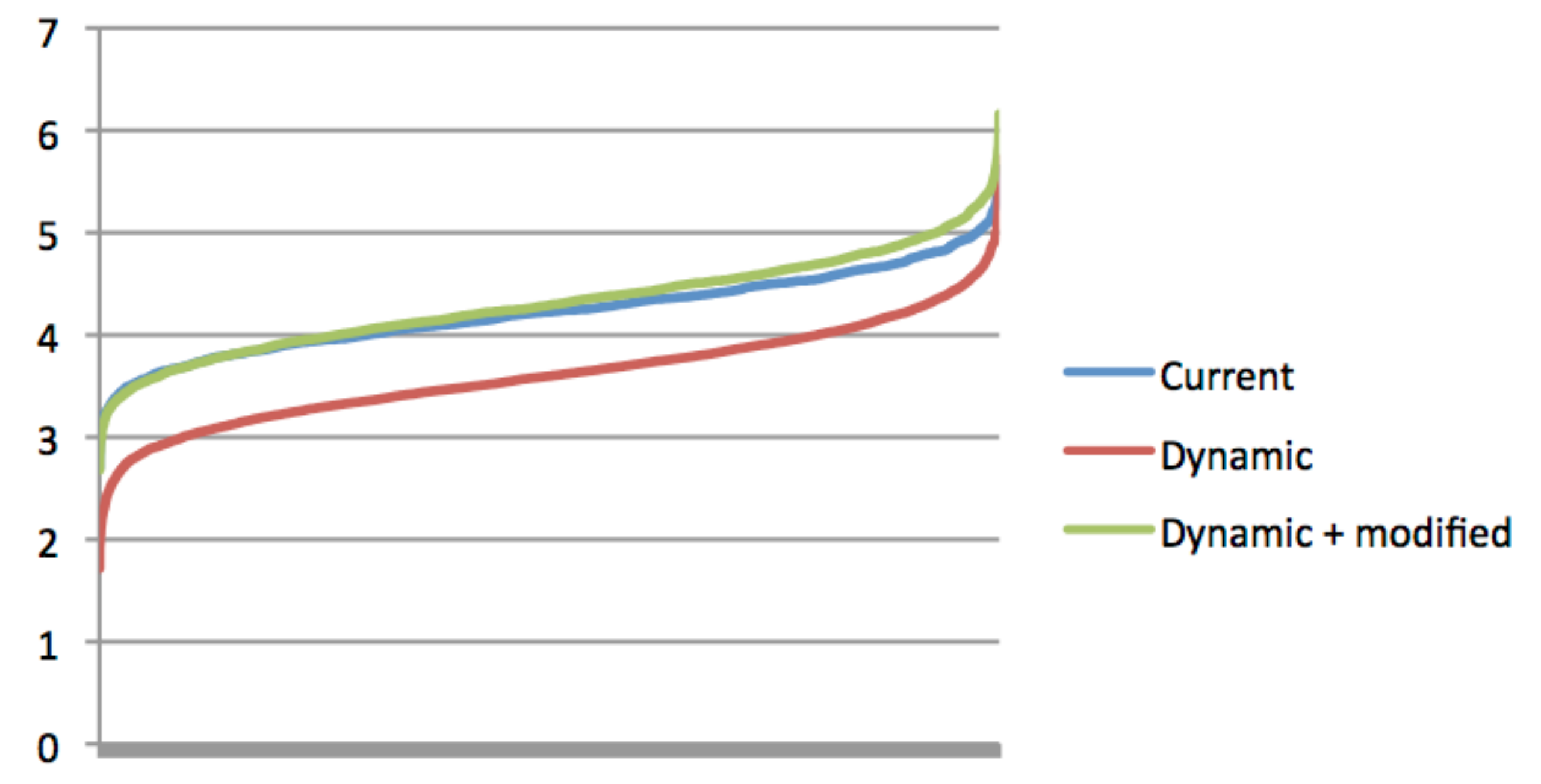}
\caption{LH measure: The dynamic method (red, average $3.61$, maximum $5.72$) 
performed better than the current Swedish method (blue, average $4.23$, maximum $5.62$ 
according to the LH measure for disproportionality between constituencies.
Included in green is for comparison the dynamic method using the modified \web.}\label{F:errorLH}
\end{figure}

\begin{figure}[htb]
\includegraphics[height=6cm]{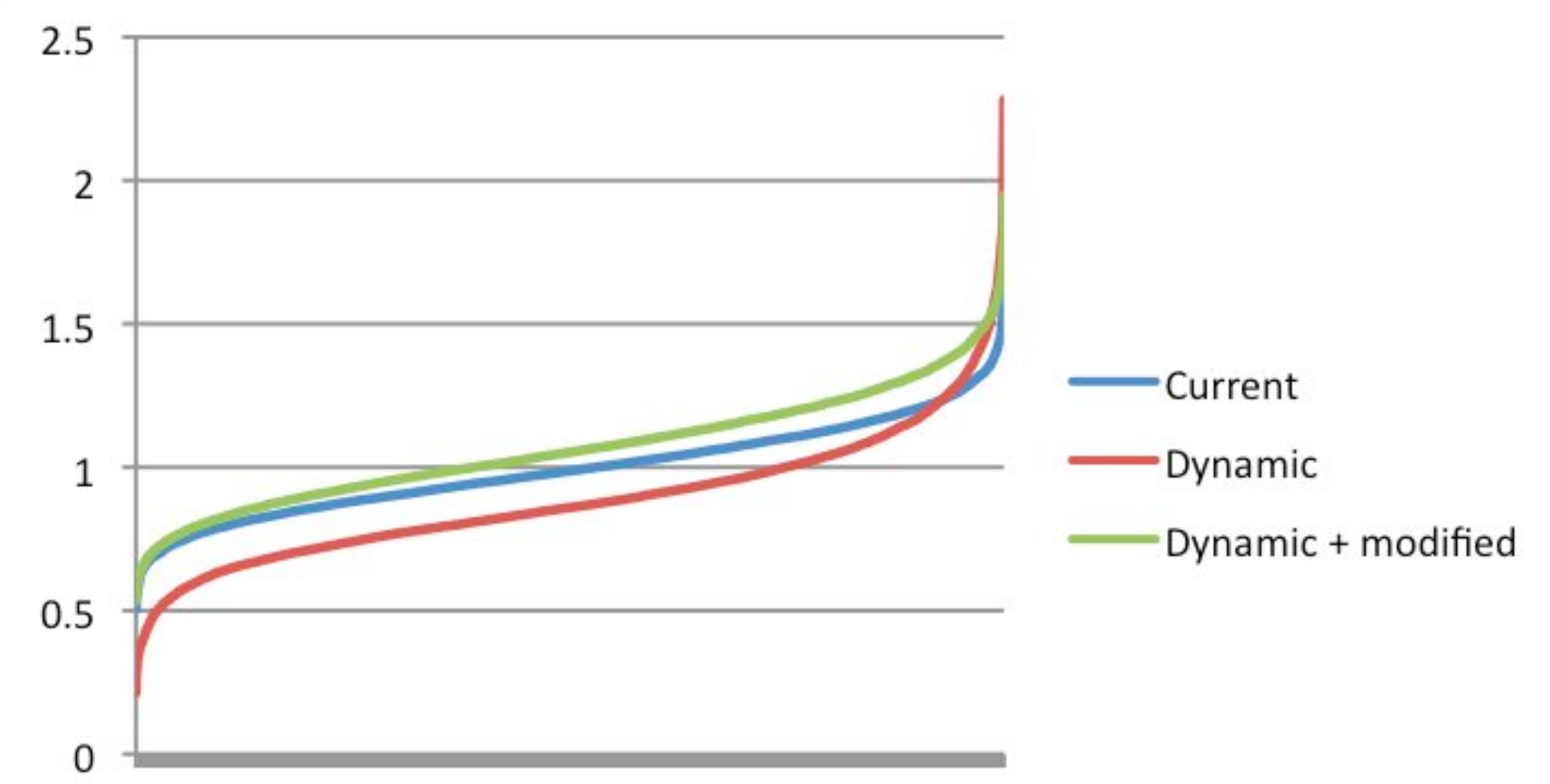}
\caption{SL measure: In the best 93\% of its performances,
the dynamic method (red, average $0.89$, maximum $2.21$) 
performed better than the current Swedish method (blue, average $0.99$, maximum $1.67$ 
according to the SL measure for disproportionality between constituencies, 
but in its worst 7\% it did worse.
Included in green is for comparison the dynamic method using the modified \web.}\label{F:errorSL}
\end{figure}

The dynamic method is constructed to give a good proportionality between parties. As the adjustment seat step is not specifically designed to enforce constituency proportionality, we wanted to know how
well the dynamic method would do with the proportionality between constituencies, in particular when the method 
generates more than the current 39 adjustment seats. 
For each of the 
simulations above we computed the constitual disproportionality with both the LH and the SL measures.
In Figures \ref{F:errorLH} and  \ref{F:errorSL} we have plotted the LH and SL measures, respectively, on the $y$-axis and 
for each of the three electoral systems separately ordered the 10 000 simulations on the $x$-axis 
by increasing measure for comparison.
As can be seen, the 
dynamic method did better than the current Swedish system.
One reason for this is the fact that we in the dynamic method use the standard \web, i.e. without the 
modifying $1.4$ factor, see the discussion at the end of this section. For comparisons sake we have also done a version of the 
dynamic method with the modifying $1.4$ factor, even though we do not suggest this as a method to be used in practice. As can 
be seen, that does increase the constituency disproportionality.
The main conclusion is that the difference in disproportionality for constituencies is not very large between the methods when 
simulating around the 2010 election result. From the perturbation results, we see that the behaviour is very much the same for the simulated election results as for the actual 
2010 voting result; in the sense that the dynamic method performs well but the current method fails to deliver party
proportionally.

%%%%%%%%%%%%%%%%%%%%%%%%%%%%%%%%%%%%%%%%%%%%%%%%%%

The disproportionality for the current Swedish system is thus on average larger than for the dynamic system. In SL-measures the 
dynamic system had however the worst 7\% cases which should be investigated further. With the LH-measure the dynamic 
method clearly outperforms the current Swedish system when it comes to proportionality between constituencies, even
though the largest three values (of 10 000) of the LH measure come from the dynamic method.

There is a reason why the current Swedish system does not behave so well for constituency proportionality.
Consider a constituency A where some party is very close to getting a permanent seat, but does not get it because of the $1.4$ factor. When the 
adjustment seats are distributed within that party, the constituency A will then be very competitive since it has so many votes that it actually should (if the
standard \web{} had been used) have gotten a permanent seat there. 
This non- symmetry, that the $1.4$ factor is used for permanent seats but not for 
adjustment seats, we believe to be the reason why we get on average improved proportionality also for constituencies with the dynamic method. 
By looking at data from old elections, it seems that the current system is somewhat favoring medium sized constituencies. We have not done any statistical tests to verify this.

\subsection{Probability of non-monotonicity for candidates}\label{S:SimNon-monotone}
In Section \ref{S:Non-monotone} we discussed the possibility of non-monotonicity for a single candidate with the dynamic method. In order to estimate the likelihood of this we checked in the simulations if this could happen.

We considered the following scenario. Let $A$ be the party that got the last overall seat in step (3), see Section \ref{S:Implement},
 and $B$ (different from $A$) the party that is closest to gaining one seat, i.e. has the highest comparison number after the last 
 seat.  If $A$ received all its seats in step (4), i.e. before the stop of permanent seats occurred then if $A$ lost that last seat to $B$ 
 the stop for permanent seats would happen earlier. If $B$ got any permanent seat in between the last two permanent seats for
  $A$ we let $K$ be the constituency in which $B$ got this seat (could be several options in which case we took the last). We ran 
  10 000 simulations as above and we got 595 triples $A,B,K$ to investigate further. 

For each such triple, we need to add votes to $B$ so that $B$ gets one more seat in total from $A$. The votes to $B$ could be 
distributed in many ways among the constituencies, which influences how likely $B$ is to get an adjustment seat in $K$. We did 
two different approaches.  First, we added all the votes in constituency $K$.
In $352$ of these cases $B$ lost a permanent seat in $K$. The 
reason they did not always lose a permanent seat is that the extra votes can cause $B$ to get the permanent seat earlier in $L$.
In 332 of the 352 cases $B$ got an extra adjustment seat in $K$ instead. 
Thus, in $20$ of the 10 000 simulations did in fact the candidate for party $B$ in $K$ lose the seat with this approach.

Secondly we tried to add the votes to $B$ more proportionally among the constituencies. To be more precise, we added votes distributed as $B$s previous votes and then one extra in constituency $K$. In this setting, the candidate for party $B$ in $K$ lost
the seat in 43 of the 10 000 simulations.

{}From both of these experiments we conclude that the non-monotonicity for a candidate is not completely esoteric. 
For election results similar to the latest Swedish results, it is sometimes possible with rather small addition of votes to create 
the non-monotonicity. When looking at the random addition of votes in one party and constituency, the non-monotonicity is 
very unlikely. First, several conditions must be fulfilled as for $A,B,K$ above and second the party $B$ must gain just the right amount 
of extra votes, properly distributed.

\section{Conclusions} \label{S:Con}
Our main conclusion of this study is that the dynamic method is an interesting alternative that should be 
seriously considered when Sweden wants to change the current system. It would have worked well in previous 
Swedish national elections and simulations around the outcome of the national election of 2010 also shows it 
would work well. It would in most cases not only give a better proportionality between parties, but also
give a better proportionality between the constituencies. This is mainly due to the fact that we discard of the 
modifying $1.4$-factor in the \web.

We have not been able to find any large drawbacks with the method. The largest is the theoretical possibility
of non-monotonicity for a given candidate for a party, as explained in \refS{S:Non-monotone} and the fact that the dynamic method would be slightly more complicated than the current method to explain to the general public. 
If desired as a safeguard against a too unfair distribution among constituencies, a practical 
implementation of the dynamic method could easily be using a lower bound on the number of permanent seats, 
either in total or for each constituency.

The simulations in \refS{S:Simulate} around the national election result of 2010 show that it was not an 
unlucky accident that the current Swedish system gave non-proportionality for parties. It does so in $95\%$ of 
the simulations, see Figure \ref{F:adjmod}. One 
conclusion is that if Sweden does not want this to happen again, the system should be changed in 
one way or the other.


\begin{thebibliography}{ABX}

\bibitem{BP} Balinski, Michel and Pukelsheim, Frederich, 'Matrices and politics'. {\em Festschrift for Tarmo Pukkila on his 60th Birthday} (Eds. E. Liski, J. Isotalo, S. Puntanen, G.P.H. Styan), Department of Mathematics, Statistics, and Philosophy: University of Tampere 2006, 233-242.  Available at {\tt http://www.math.uni-augsburg.de/stochastik/pukelsheim/2006d.pdf}

\bibitem{BY} Balinski, Michel and Young, Peyton, {\it Fair Representation, 2nd ed}, Brookings Institution (2001).

\bibitem{BF} Brams, S.J., Fishburn, P.C., 'Proportional representation in Variable-Size Legislatures', {\it Social choice and wellfare}, {\bf 1}, (1984), 211--229.

\bibitem{Fr} Fr\"oberg and Sundstr\"om, appendix in {\em Partiell f\"orfattningsreform}, SOU:1967:26.

\bibitem{G91} Gallagher, Michael, 'Proportionality, disproportionality and Electoral Systems', 
{\it Electoral Studies}, {\bf 10}, no. 1, 33--51. (1991).

\bibitem{Lanke} Lanke, Jan, 'What is it that makes the Swedish election act fail?', talk at workshop on electoral methods
at KTH May 30-31, 2011. See http://www.math.kth.se/wem/Lanke.pptx.

\bibitem{JL} Janson, Svante and Linusson, Svante, 'The probability of the Alabama paradox', 
{\em J. of Applied Prob.}, {\bf 49}, No. 3, pp. 773-794 (2012). %Arxiv:math/1104.2137.

\bibitem{LH}Loosemore, J. and Hanby, V. (1971): 'The Theoretical Limits of Maximum Distortion: Some Analytic Expressions for Electoral Systems': {\em British Journal of Political Science} 1, 467--477.

\bibitem{Ni} Nilsson, Bengt-{\AA}ke, et al, {\em Proportionalitet i val samt f{\"o}rhandsanm{\"a}lan av partier och kandidater}, SOU 2012:94 (2013). {\tt http://www.regeringen.se/sb/d/15631/a/206860}

\bibitem{EU} Office for Promotion 
of Parliamentary Democracy {\em Electoral systems - The link between 
governance, elected members and voters} , European Parliament (2011).
{\tt http://www.europarl.europa.eu/pdf/oppd/Page\_8/Electoral-systems-LR-for-WEB.pdf}

\bibitem{Sydow}von Sydow, Bj\"orn, {\em V\"agen till enkammarriksdagen}, Tiden (1989).

\bibitem{Val} Valmyndigheten, Swedish Election system in brief,  {\tt 
http://www.val.se/pdf/electionsinsweden\_webb.pdf}.

\end{thebibliography}
\end{document}